\documentclass[journal]{IEEEtran}
\IEEEoverridecommandlockouts
\usepackage{amsmath,amssymb,latexsym,epsfig,psfrag,graphicx,cases}
\newtheorem{prp}{Proposition}
\newtheorem{thm}[prp]{Theorem}
\newtheorem{cor}[prp]{Corollary}
\newtheorem{lem}[prp]{Lemma}

\newenvironment{pf}{\noindent{\it Proof. }}{\hfill$\Box$\smallbreak}
\newenvironment{pf*}[1]{\smallbreak\noindent{\it #1}}{\hfill$\Box$\smallbreak}
\newcounter{remark}
\newenvironment{rmk}{\addtocounter{remark}{1}\smallbreak\noindent
  {\em Remark \theremark.}}{\smallbreak}
\newcommand{\realR}{\mathbb{R}}

\begin{document}
\title{Linear Quadratic Dual Control}
  \author{Anders Rantzer
  \thanks{The author is affiliated with Automatic Control LTH, Lund
    University, Box 118, SE-221 00 Lund, Sweden. He is a member of the Excellence Center ELLIIT and Wallenberg AI, Autonomous Systems and Software Program (WASP). Support was received from the European Research Council (Advanced Grant 834142) }}
\maketitle

\begin{abstract}%
This is a draft paper originally posted on Arxiv as a documentation of a plenary lecture at CDC2023. The core material has been accepted for publication at L4DC 2024. 

Certainty equivalence adaptive controllers are analysed using a ``data-driven Riccati equation'', corresponding to the model-free Bellman equation used in Q-learning. The equation depends quadratically on data correlation matrices. This makes it possible to derive simple sufficient conditions for stability and robustness to unmodeled dynamics in adaptive systems. The paper is concluded by short remarks on how the bounds can be used to quantify the interplay between excitation levels and robustness.
\end{abstract}

\section{Introduction}
The history of adaptive control dates back at least to aircraft autopilot development in the 1950s. Following the landmark paper \cite{ast+wit73}, a surge of research activity during the 1970s derived conditions for convergence, stability, robustness and performance under various assumptions. For example, \cite{ljung77ana} analysed adaptive algorithms using averaging, \cite{goodwin1981discrete} derived an algorithm that gives mean square stability with probability one, while \cite{guo1995convergence} gave conditions for the optimal asymptotic rate of convergence. On the other hand, conditions that may cause instability were studied in \cite{ega79book}, \cite{ioannou1984instability} and \cite{rohrs1985robustness}. Altogether, the subject has a rich history documented in numerous textbooks, such as \cite{aastrom2013adaptive}, \cite{goodwin2014adaptive}, and \cite{sastry2011adaptive}. 

Recently, there has been a renewed interest in analysis of adaptive controllers, driven by progress in statistical machine learning. See \cite{tsiamis2023statistical} for a review. In parallel, there is also a rapidly developing literature on (off-line) data-driven control. \cite{de2019formulas,markovsky2021behavioral,berberich2020data}.

In this paper, the focus is on worst-case models for disturbances and uncertain parameters, as discussed in \cite{cusumano1988nonlinear,sun1987theory,vinnicombe2004examples,2004megretskinonlinear} and more recently in \cite{rantzer2021minimax,cederberg2022synthesis,kjellqvist2022minimax}. However, the disturbances in this paper are assumend to be bounded in terms of past states and inputs. This causality constraint is different from above mentioned references.

\section{Notation}
The set of $n\times m$ matrices with real coefficients is denoted $\realR^{n\times m}$. The transpose of a matrix $A$ is denoted $A^\top$. For a symmetric matrix $A\in\realR^{n\times n}$, we write $A\succ0$ to say that $A$ is positive definite, while $A\succeq0$ means positive semi-definite. 
Given $x\in\realR^n$ and $A\in\realR^{n\times n}$, the notation $|x|^2_A$ means 
$x^\top Ax$. The expression $\min_K
  {\tiny\begin{bmatrix}I\\K\end{bmatrix}^\top}Q
  {\tiny\begin{bmatrix}I\\K\end{bmatrix}}$ is equivalent to $Q^{xx}-Q^{xu}(Q^{uu})^{-1}Q^{ux}$ when $Q={\tiny\begin{bmatrix}Q^{xx}&Q^{xu}\\Q^{ux}&Q^{uu}\end{bmatrix}}\succ0$.

\section{A Data-driven Riccati Equation}
Consider a linear quadratic optimal control problem:
\begin{align*}
  \text{Minimize }&\sum_{t=0}^\infty\left(|x_t|^2+|u_t|^2\right)\\
  \text{subject to }&x_{t+1}=Ax_t+Bu_t, 
\end{align*}
with $x_0\in\realR^n$ given and $u_t\in\realR^m$. Assuming that the system is stabilizable, the optimal value has the form $|x_0|^2_{P}$ where $P\in\realR^{n\times n}$ can be obtained by solving the Riccati equation
\begin{align}
  P&=\min_{K}\left[I+K^\top K+(A+BK)^\top P(A+BK)\right].
\label{eqn:RicP}
\end{align}
Define $Q:=I+\begin{bmatrix}A&B\end{bmatrix}^\top P\begin{bmatrix}A&B\end{bmatrix}\in\realR^{(n+m)\times(n+m)}$ . Then \eqref{eqn:RicP} can be written as
\begin{align}
  Q-I= \begin{bmatrix}A&B\end{bmatrix}^\top\min_{K}\left(\begin{bmatrix}I\\K\end{bmatrix}^\top Q\begin{bmatrix}I\\K\end{bmatrix}\right)\begin{bmatrix}A&B\end{bmatrix}
\label{eqn:Ric}
\end{align}
or alternatively
\begin{align}
  {\begin{bmatrix}x\\u\end{bmatrix}^\top (Q-I)
  \begin{bmatrix}x\\u\end{bmatrix}}
  &=x_+^\top\min_{K}\left(\begin{bmatrix}I\\K\end{bmatrix}^\top Q\begin{bmatrix}I\\K\end{bmatrix}\right)x_+
\label{eqn:RicQ}
\end{align}
where $x_+=Ax+Bu$. The equation \eqref{eqn:RicQ} is sometimes called \emph{model free}, since it does not include the model parameters $(A,B)$. This makes it possible to collect data points $(x,u,x_+)$ and use \eqref{eqn:RicQ} to get information about $Q$. In fact, the total matrix $Q$ can be computed from a trajectory $x_0,u_0,\ldots,x_{t-1},u_{t-1},x_t$ spanning all directions of $(x_t,u_t)$, using the equation
\begin{align*}
  &\begin{bmatrix}x_0&\ldots&x_{t-1}\\u_0&\ldots&u_{t-1}\end{bmatrix}^\top(Q-I)
  \begin{bmatrix}x_0&\ldots&x_{t-1}\\u_0&\ldots&u_{t-1}\end{bmatrix}\\
  &=
  \begin{bmatrix}x_1&\ldots&x_t\end{bmatrix}^\top
  \left(\begin{bmatrix}I\\K\end{bmatrix}^\top Q\begin{bmatrix}I\\K\end{bmatrix}\right)
  \begin{bmatrix}x_1&\ldots&x_t\end{bmatrix}
\end{align*}
This is essentially equation (3) in \cite{bradtke1992reinforcement} and (14) in \cite{rizvi2018output}. The fact that the equation is over-determined for $t>n+m$  makes it natural to multiply from the left and right by
\begin{align*}
  \begin{bmatrix}\lambda^tx_0&\lambda^{t-1}x_1&\ldots&x_{t-1}\\\lambda^tu_0&\lambda^{t-1}u_1&\ldots&u_{t-1}\end{bmatrix},
\end{align*}
where $\lambda\in(0,1)$ is a forgetting factor. With notation 
\begin{align}
  {\Sigma}_t&=\begin{bmatrix}
        \Sigma_t^{xx}&\Sigma_t^{xu}\\
        \Sigma_t^{ux}&\Sigma_t^{uu}
      \end{bmatrix}=\sum_{k=0}^{t-1}\lambda^{t-1-k} \begin{bmatrix}x_k\\u_k\end{bmatrix}
  \begin{bmatrix}x_k\\u_k\end{bmatrix}^\top+\lambda^t\Sigma_0\label{eqn:Sigt}\\
  \hat{{\Sigma}}_t&=\sum_{k=0}^{t-1}\lambda^{t-1-k}{x_{k+1}
    \begin{bmatrix}x_k\\u_k\end{bmatrix}^\top,}\label{eqn:hatSigt}
\end{align} 
where we have added a regularization term $\lambda^t\Sigma_0\succ0$ to make $\Sigma_t$ intervible, 
this gives a \emph{data driven Riccati equation}
\begin{align}
 {\Sigma}_t\left(Q_t-I\right){\Sigma}_t
 &=\hat{{\Sigma}}_t^\top\min_{K}\left(\begin{bmatrix}I\\K\end{bmatrix}^\top Q_t\begin{bmatrix}I\\K\end{bmatrix}\right)\hat{{\Sigma}}_t.
\label{eqn:Bellmantau}
\end{align}
The equation \eqref{eqn:Bellmantau} makes it possible to estimate the optimal cost matrix $Q$ directly from data. The model-free aspect should not be over-emphasized, since an algebraically equivalent approach is to first estimate $A$ and $B$ according to $\begin{bmatrix}\hat{A}_t&\hat{B}_t\end{bmatrix}:=\hat{{\Sigma}}_t{\Sigma}_t^{-1}$ and then solve the standard Riccati equation
\begin{align*}
  P_t&=\min_{K}\left[I+K^\top K+(\hat{A}_t+\hat{B}_tK)^\top P_t(\hat{A}_t+\hat{B}_tK)\right]
\end{align*}
in analogy with \eqref{eqn:RicP}.
Nevertheless, the fact that \eqref{eqn:Bellmantau} is quadratic in the correlation matrices $\hat{{\Sigma}}_t$, ${\Sigma}_t$ makes it a good starting point for analysis of data-driven controllers.

\section{Problem formulation}
\label{sec:controller}
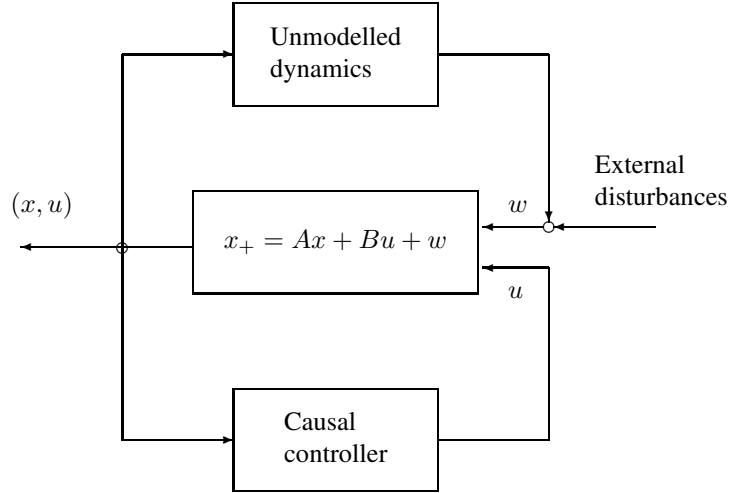
\begin{figure}
  \begin{center}
  \setlength{\unitlength}{.009mm}%
\begin{picture}(9441,7266)(1318,-8494)
\put(8551,-5611){$u$}
\put(2851,-4861){\circle{150}}
\put(4501,-2761){\framebox(3000,1500){\begin{tabular}{ll}
        Unmodelled\\ dynamics
      \end{tabular}}}
\put(4501,-8461){\framebox(3000,1500){\begin{tabular}{ll}
        Causal\\ controller
      \end{tabular}}}
\put(3901,-5536){\framebox(4200,1500){$x_+=Ax+Bu+w$}}
\put(7501,-7711){\line( 1, 0){1650}}
\put(9151,-7711){\line( 0, 1){2550}}
\put(9151,-5161){\vector(-1, 0){975}}
\put(7501,-2011){\line( 1, 0){1650}}
\put(9151,-2011){\vector( 0,-1){2475}}
\put(9076,-4561){\vector(-1, 0){900}}
\put(2851,-4861){\line(1, 0){1050}}
\put(2851,-4861){\line( 0, 1){2850}}
\put(2851,-2011){\vector( 1, 0){1650}}
\put(2851,-4861){\vector(-1, 0){1500}}
\put(2851,-4861){\line( 0,-1){2850}}
\put(2851,-7711){\vector( 1, 0){1650}}
\put(10726,-4561){\vector(-1, 0){1500}}
\put(1201,-4336){$(x,u)$}
\put(8551,-4336){$w$}
\put(9600,-4000){\begin{tabular}{ll}
        External\\ disturbances
      \end{tabular}}
\put(9151,-4561){\circle{150}}
\end{picture}%
  \end{center}
  \caption{The system \eqref{eqn:state} is connected to a causal controller defined by \eqref{eqn:controller}. The disturbance $w$ could represent a combination of external disturbances and unmodeled dynamics. }
\label{fig:blockdiag}
\end{figure}

For the linear system 
\begin{align}
  x_{t+1}&=Ax_t+Bu_t+w_t&t&\ge0\label{eqn:state}
\end{align}
we will analyse controllers of the following form (See Figure~\ref{fig:blockdiag}):
\begin{align}
  \begin{cases}
    \Sigma_t=\lambda\Sigma_{t-1}+\begin{bmatrix}x_{t-1}^\top&u_{t-1}^\top\end{bmatrix}^\top\begin{bmatrix}x_{t-1}^\top&u_{t-1}^\top\end{bmatrix},
    \qquad \Sigma_0\succ0\\
    \hat{\Sigma}_t=\lambda\hat{\Sigma}_{t-1}+x_t\begin{bmatrix}x_{t-1}^\top&u_{t-1}^\top\end{bmatrix},
    \qquad\qquad\qquad\quad\;\, \hat\Sigma_0=0\\
    \, u_t=K_tx_t+\epsilon_t
  \end{cases}
\label{eqn:controller}
\end{align}
where $\lambda\in(0,1)$ and $K_t$ is obtained from $\Sigma_t$ and $\hat{\Sigma}_t$ as the minimizing argument in a solution of \eqref{eqn:Bellmantau} and $\epsilon_t$ is added to provide excitation (exploration). In the adaptive control literature such controllers are known as \emph{certainty equivalence} controllers, since the main feedback term $K_tx_t$ can be obtained by using the model estimate $\begin{bmatrix}\hat{A}_t&\hat{B}_t\end{bmatrix}:=\hat{{\Sigma}}_t{\Sigma}_t^{-1}$ for optimization, disregarding its uncertainty.

We will give conditions that guarantee stability and performance (in terms of $\sum_t(|x_t|^2+|u_t|^2)$) for the closed loop system.
Obviously, this cannot be done without some prior constraints on the model (such as stabilizability). Many alternative formulations have been considered in the literature,  \cite{konstantinov1993perturbation},\cite[Assumption~2]{mania2019certainty}, but the following simple version will be enough for our purposes. \emph{Given any $\beta>1$, we restrict the model parameters to be in the set $\mathcal{M}_{\beta}$ consisting of all pairs $(A,B)$ such that the Riccati equation \eqref{eqn:Ric}
has a solution $Q$ with $I\preceq Q\preceq\beta^2I$.}

We also need to put assumptions on the disturbance $w$. For this purpose, introduce the notation
\begin{align}
  \begin{bmatrix}\Sigma_t^{wx}&\Sigma_t^{wu}\end{bmatrix}
  &:=\sum_{k=0}^{t-1}\lambda^{t-1-k}w_k\begin{bmatrix}x_k^\top&u_k^\top\end{bmatrix}-\lambda^t\begin{bmatrix}A&B\end{bmatrix}\Sigma_0
  \label{eqn:tildeSigt}
\end{align} 
and notice that $\begin{bmatrix}\Sigma_t^{wx}&\Sigma_t^{wu}\end{bmatrix}=\hat{\Sigma}_t-\begin{bmatrix}A&B\end{bmatrix}\Sigma_t$.
In statistical system analysis, $w_t$ is often modeled as white noise with fixed amplitude. In this paper, we will instead think of $w_t$ as unmodeled dynamics, which makes it natural to bound $\Sigma_t^{wx}$ and $\Sigma_t^{wu}$ in terms of $\Sigma_t^{xx}$ and $\Sigma_t^{uu}$.

\section{Main Results}

Our first and main theorem is quantifying the effect in a single time-step of using the feedback gain $K_t$ obtained from the data-driven Riccati equation instead of the optimal gain from the standard Riccati equation:

\begin{thm}
  Consider $\beta,\rho\in\realR_+$ with $2\beta^2\rho(\rho+2)<1$, a matrix $\hat{{\Sigma}}_t$ and a positive definite matrix $\Sigma_t$.
  Let $(A,B)\in\mathcal{M}_\beta$ with $P$ satisfying Riccati equation \eqref{eqn:RicP}, while
  \begin{align}
\left\|\begin{bmatrix}A&B\end{bmatrix}-\hat{\Sigma}_t\Sigma_t^{-1}\right\|&\le\rho.
  \label{eqn:preceq}
  \end{align} 
  Then \eqref{eqn:Bellmantau}
  has a solution with minimizer $K_t$ and
  {\begin{align}
   \frac{1}{1-2\beta^2\rho(\rho+2)}P&\succeq I+K_t^\top K_t+
   (A+BK_t)^\top P(A+BK_t).
   \label{eqn:main3}
   \end{align}}
\label{thm:dataRic}
\end{thm}


\begin{rmk}
In analysis of \eqref{eqn:state}, the assumption \eqref{eqn:preceq} can be written as $\big\|\begin{bmatrix}\Sigma_t^{wx}&\Sigma_t^{wu}\end{bmatrix}\Sigma_t^{-1}\big\|\le\rho$. This inequality tends to fail for small values of $t$, but proper choice of excitation and regularization together with sufficiently small disturbances makes it valid after collection of enough data.
\end{rmk}

\begin{rmk}
Notice that the optimal control law $\bar{K}$ from the standard Riccati equation satisfies 
\begin{align}
  P&=I+\bar{K}^\top \bar{K}+(A+B\bar{K})^\top P(A+B\bar{K})
\label{eqn:Lyap}
\end{align}
while $K_t$ from the data driven Riccati equation gives
\begin{align*}
  P_t&=I+K_t^\top K_t+(\hat{A}_t+\hat{B}_tK_t)^\top P_t(\hat{A}_t+\hat{B}_tK_t)
\end{align*}
The purpose of Theorem~\ref{thm:dataRic} is to show what happens to the storage function $x^\top Px$ when we use $K_t$ instead of the optimal $\bar{K}$ to control the system. As expected, the difference vanishes as \hbox{$\rho\to0$}. Larger values of $\beta$ open up for systems that are harder to stabilizable, which means that $\rho$ needs to be smaller and more accurate data is needed. This is reflected in the left hand side of \eqref{eqn:main3}.
\end{rmk}

  The optimal closed loop system obtained with $u_t=\bar Kx_t$ is stable, with the Lyapunov function $|x_t|^2_P$ decays according to 
  $|x_{t+1}|^2_P-|x_t|^2_P\le -|x_t|^2-|\bar Kx_t|^2$. The strict decay rate gives robustness to deviations from optimality. Hence, as long as the difference between the left hand sides of \eqref{eqn:main3} and \eqref{eqn:Lyap} is not larger than the identity matrix, stability is retained. This is exploited in to get the following corollary.

\begin{cor}
  Consider $(A,B)\in\mathcal{M}_\beta$ with $P$ being a positive definite solution to the Riccati equation \eqref{eqn:RicP}, while $\beta,\rho, \gamma$ are positive numbers. Let
  \begin{align*}
    x_{t+1}&=Ax_t+Bu_t+w_t,&
    u_t&=K_tx_t+\epsilon_t
  \end{align*}
  and suppose that $\Sigma_t,\hat{\Sigma}_t$ are defined by \eqref{eqn:Sigt}-\eqref{eqn:hatSigt}, while \eqref{eqn:preceq} holds
  for $t\ge t_0$. For each $t$, let $Q_t$ be a solution of \eqref{eqn:Bellmantau} with $K_t$ being a minimizing argument. Define 
  \begin{align*}
    \alpha&:=\beta^2+\frac{1}{1-\beta^2/\gamma^2}\left(1-\frac{\beta^2}{1-2\beta^2\rho(\rho+2)}\right).
  \end{align*}
  If $\alpha>0$, then
  \begin{align}
    \sum_{t=t_0}^{T-1}\left(|x_t|^2+|K_tx_t|^2\right)
    &\le\alpha^{-1}|x_{t_0}|^2_P+\frac{\gamma^2}{\alpha}\sum_{t=t_0}^{T-1}|B\epsilon_t+w_t|^2.
  \label{eqn:gainbound}
  \end{align}
\label{thm:main}
\end{cor}

\section{Concluding remarks}
The results above are the first steps towards a proof that a properly designed certainty equivalence controller can robustly stabilize the set $\mathcal{M}_\beta$ for arbitrarily large values of $\beta$. Here robustness means tolerance to a certain amount of unmodeled dynamics. The argument, with further details given in \cite{rantzer2023linear}, has the following steps:
\begin{enumerate}
  \item Pick any $\beta>1$.
  \item Select $\rho>0$ small enough to make $[1-2\beta^2\rho(\rho+2)]^{-1}<1+\beta^{-2}$. Then Corollary~\ref{thm:main} proves finite gain from $B\epsilon+w$ to $x$, provided that $\big\|\begin{bmatrix}\Sigma_t^{wx}&\Sigma_t^{wu}\end{bmatrix}\Sigma_t^{-1}\big\|\le\rho$.
  \item Design the excitation term $\epsilon_t$ to bound the impact of $\Sigma_t^{-1}$. 
  \item Introduce a bound on $w_t$ to limit the size of $\begin{bmatrix}\Sigma_t^{wx}&\Sigma_t^{wu}\end{bmatrix}$.
\end{enumerate}

\section{Acknowledgement}The author is a member of the Excellence Center ELLIIT and Wallenberg AI, Autonomous Systems and Software Program (WASP). Support was received from the European Research Council (Advanced Grant 834142)





\section{Appendix: Proofs}
\begin{lem}
  Consider positive numbers $\beta$ and $\rho$, matrices $\hat\Sigma$, $\tilde\Sigma$ and positive semi-definite matrices $\Sigma$, $P$, $Q$. Suppose that $(\hat\Sigma-\tilde\Sigma)^\top P(\hat\Sigma-\tilde\Sigma)=\Sigma(Q-I)\Sigma$ and  $\tilde\Sigma^\top\tilde\Sigma\le\rho^2\Sigma^2$. If $I\preceq Q\preceq\beta^2I$, then $\hat\Sigma^\top P\hat\Sigma\preceq \Sigma Q\Sigma+\left(\beta^2\rho(\rho+2)-1\right)\Sigma^2$.
\label{lem:Ricpert0}
\end{lem}

\begin{pf}
  For every $z$ with appropriate dimension 
  \begin{align*}
    |\hat\Sigma z|^2_P&\le\left(|(\hat\Sigma-\tilde\Sigma)z|_P+|\tilde\Sigma z|_P\right)^2\\
    &=\left(|\Sigma z|_{Q-I}+|\tilde\Sigma z|_P\right)^2\\
    &\le\left(|\Sigma z|_{Q-I}+\beta\rho|\Sigma z|\right)^2\\
    &=|\Sigma z|^2_{Q-I}+\beta^2\rho^2|\Sigma z|^2+2|\Sigma z|_{Q-I}\beta\rho|\Sigma z|\\
    &\le|\Sigma z|^2_{Q-I}+\beta^2\rho^2|\Sigma z|^2+2\beta^2\rho|\Sigma z|^2\\
    &=|\Sigma z|^2_{Q}+\left(\beta^2\rho(\rho+2)-1\right)|\Sigma z|^2.
  \end{align*}
  Here the first inequality is the triangle inequality, the second follows as $P\preceq\beta^2I$ and $\tilde\Sigma^\top\tilde\Sigma\preceq\rho^2\Sigma^2$, while the third inequality follows from the assumption $I\preceq Q\preceq\beta^2I$. This completes the proof.
\end{pf}

\begin{lem}
  Consider matrices $A,B,\bar K$ and $\bar Q\succeq I$.
  Suppose 
  \begin{align}
    \begin{bmatrix}A&B\end{bmatrix}^\top \begin{bmatrix}I\\\bar K\end{bmatrix}^\top \bar Q\begin{bmatrix}I\\\bar K\end{bmatrix}\begin{bmatrix}A&B\end{bmatrix}&\preceq \bar Q-I.
  \label{eqn:Sigineq}
  \end{align}
  Then there exists a unique $Q$ with $I\preceq Q\preceq\bar{Q}$ such that
  \begin{align}
    \begin{bmatrix}A&B\end{bmatrix}^\top \min_{K}\bigg(\begin{bmatrix}I\\K\end{bmatrix}^\top Q\begin{bmatrix}I\\K\end{bmatrix}\bigg)\begin{bmatrix}A&B\end{bmatrix}&= Q-I.
  \label{eqn:Qric}
  \end{align}
\label{lem:Ricsol}
\end{lem}

\begin{pf}
  The matrix
  \begin{align*}
    \bar P&:=\begin{bmatrix}I&\bar{K}^\top\end{bmatrix} \bar Q\begin{bmatrix}I&\bar{K}^\top\end{bmatrix}^\top
  \end{align*}
  satisfies
  \begin{align*}
    \bar P&=(A+B\bar K)^\top\bar P(A+B\bar K)+I+\bar K^\top \bar K,
  \end{align*}
  so $P_0=\bar P$ and the Riccati iteration
  \begin{align*}
    P_{k+1}&=\min_K\left[(A+B\bar K)^\top P_k(A+B\bar K)+I+\bar K^\top \bar K\right]
  \end{align*}
  gives a decreasing sequence with limit $P_\infty\succeq I$. Setting 
  \begin{align*}
    Q&:=\begin{bmatrix}A&B\end{bmatrix}^\top P_\infty\begin{bmatrix}A&B\end{bmatrix}
  \end{align*}
  proves existence. Stabilizing solutions to the Riccati equation are unique \cite[Ch.~3]{brockett2015finite}, so also $Q$ must be unique.
\end{pf}

\begin{pf*}{Proof of Theorem~\ref{thm:dataRic}.}
 By definition, the assumption $(A,B)\in\mathcal{M}$ implies existence of $Q$ satisfying \eqref{eqn:Ric} and $I\preceq Q\preceq\beta^2I$. Let $\bar{K}$ be the corresponding minimizing argument in \eqref{eqn:Ric}.
 Define
 \begin{align*}
   P&:=\begin{bmatrix}I&\bar{K}^\top\end{bmatrix} Q\begin{bmatrix}I&\bar{K}^\top\end{bmatrix}^\top\\
   \tilde{\Sigma}&=\hat{\Sigma}_t-\begin{bmatrix}A&B\end{bmatrix}\Sigma_t.
 \end{align*}
 Then \eqref{eqn:Ric} gives
 $(\hat\Sigma_t-\tilde\Sigma_t)^\top P(\hat\Sigma_t-\tilde\Sigma_t)=\Sigma_t(Q-I)\Sigma_t$,
 so Lemma~\ref{lem:Ricpert0} implies
  \begin{align*}
    \hat\Sigma_t^\top P\hat\Sigma_t&\preceq \Sigma_t(Q-\check\alpha I)\Sigma_t,
  \end{align*}
  where $\check\alpha:=1-\beta^2\rho(\rho+2)$. It follows that \eqref{eqn:Sigineq} holds with $\bar{Q}=\check\alpha^{-1}Q$, so  Lemma~\ref{lem:Ricsol} proves existence of $Q_t$ (and $K_t$) satisfying \eqref{eqn:Bellmantau} and $I\preceq Q_t\preceq\check\alpha^{-1}Q$. 
  
  Next define 
  \begin{align*}
     P_t&:=\begin{bmatrix}I&K_t^\top\end{bmatrix} Q\begin{bmatrix}I&K_t^\top\end{bmatrix}^\top\\
     \bar\Sigma_t&:=\begin{bmatrix}A&B\end{bmatrix}\Sigma_t.
  \end{align*}
  Then the inequalities $(\bar\Sigma_t+\tilde\Sigma_t)^\top P_t(\bar\Sigma_t+\tilde\Sigma_t)=\Sigma_t(Q_t-I)\Sigma_t$ and $Q_t\preceq\alpha^{-1}\beta^2I$ allow application of Lemma~\ref{lem:Ricpert0} to get
  \begin{align}
    \bar\Sigma_t^\top P_t\bar\Sigma_t= \Sigma_t(Q_t-\hat{\alpha}I)\Sigma_t,
  \label{eqn:QRic2}
  \end{align}
  with $\hat\alpha:=1-\check\alpha^{-1}\beta^2\rho(\rho+2)$. It follows from  Lemma~\ref{lem:Ricsol} that $Q\preceq\hat\alpha^{-1}Q_t$, so
  \begin{align*}
    &(A+BK_t)^\top P(A+BK_t)\\
    &=\begin{bmatrix}I\\K_t\end{bmatrix}^\top\begin{bmatrix}A&B\end{bmatrix}^\top\begin{bmatrix}I\\\bar{K}\end{bmatrix}^\top Q\begin{bmatrix}I\\\bar{K}\end{bmatrix}\begin{bmatrix}A&B\end{bmatrix}
    \begin{bmatrix}I\\K_t\end{bmatrix}\\
    &=\begin{bmatrix}I\\K_t\end{bmatrix}^\top(Q-I)\begin{bmatrix}I\\K_t\end{bmatrix}\\
    &\preceq\begin{bmatrix}I\\K_t\end{bmatrix}^\top(\hat{\alpha}^{-1}Q_t-I)\begin{bmatrix}I\\K_t\end{bmatrix}\\
    &=\hat{\alpha}^{-1}P_t-(I+K_t^\top K_t)\\
    &\preceq\hat{\alpha}^{-1}\check\alpha^{-1}P-(I+K_t^\top K_t)\\
    &=\left[1-2\beta^2\rho(\rho+2)\right]^{-1}P-(I+K_t^\top K_t)
  \end{align*}
  and the proof is complete.
  \end{pf*}

\begin{pf*}{Proof of Corollary~\ref{thm:main}.}
  \begin{align*}
    &\!\!\!\!\!\!\max_{w_t}\left(|x_{t+1}|_P^2-\gamma^2|B\epsilon_t+w_t|^2\right)\\
    &=\max_{w_t}\left(|(A+BK_t)x_t+w_t|_P^2-\gamma^2|w_t|^2\right)\\
    &=|(A+BK_t)x_t|_{P(I-\gamma^{-2}P)^{-1}}^2\\
    &\le\left(1-\beta^2/\gamma^2\right)^{-1}|(A+BK_t)x_t|_P^2\\
    &\le\frac{1}{1-\beta^2/\gamma^2}
    \left(\frac{|x_t|_P^2}{1-2\beta^2\rho(\rho+2)}-|x_t|^2-|K_tx_t|^2\right)\\
    &\le|x_t|_P^2+\left(\frac{1}{(1-\beta^2/\gamma^2)(1-2\beta^2\rho(\rho+2))}-1\right)\beta^2|x_t|^2\\
    &\quad-\frac{|x_t|^2+|K_tx_t|^2}{1-\beta^2/\gamma^2}\\
    &\le|x_t|_P^2-\alpha(|x_t|^2+|K_tx_t|^2)
  \end{align*}
  Summing over $t_0\le t\le T-1$ gives the desired result.
\end{pf*}

\end{document}